\documentclass[a4paper,12pt]{amsart}
\usepackage{amsmath}
\usepackage{amsfonts}
\usepackage{amssymb}
\usepackage{amsthm}
\newtheorem{thm}{Theorem}[section]

\numberwithin{equation}{section}

\newcommand{\Complex}{\mathbb C}

\newcommand{\To}{\rightarrow}

\newcommand{\pfbundle}[4]{#1 \rightsquigarrow #2 \xrightarrow{\smash[t]{#3}} #4}

\newcommand{\cat}{{\star}}
\newcommand{\bdelt}[1]{\mathsf{#1}}

\author{Kim Morrison}
\title{An Introduction to Pull-backs of Bundles and Homotopy Invariance}
\date{August 26 2000}

\begin{document}

\begin{abstract}
We define the pull-back of a smooth principal fibre bundle, and
show that it has a natural principal fibre bundle structure. Next,
we analyse the relationship between pull-backs by homotopy
equivalent maps. The main result of this article is to show that
for a principal fibre bundle over a paracompact manifold, there is
a principal fibre bundle isomorphism between pull-backs obtained
from homotopic maps. This enables simple proofs of several results
on the structure of principal fibre bundles. No new results are
obtained---this is simply an accessible presentation of an
important idea.
\end{abstract}

\maketitle

\section{The pull-back}
In this section we give a brief definition of the pull-back of a
principal fibre bundle \cite{cho:amp,hus:fb}. For brevity many
details are left to the references \cite{ish:mdgfp}. Suppose $\xi
= \pfbundle{G}{P}{\pi}{N}$ is a principal fibre bundle with
structure group $G$ over the base space $N$. If $f:M \To N$, we
define the \emph{pull-back of $\xi$ by $f$} as
\begin{eqnarray*}
    E & = & \left\{ (x,\bdelt{p}) \in M \times P \mid f(x)=\pi(\bdelt{p}) \right\}, \\
    & & \pi_E : E \To M \quad \textrm{by} \quad \pi_E(x,\bdelt{p})=x, \\
    f^*\xi & = & \pfbundle{G}{E}{\pi_E}{M}.
\end{eqnarray*}
The right action of $G$ on $E$ is given simply by
$(x,\bdelt{p})g=(x,\bdelt{p} g)$. The fibre bundle structure is
also simple to describe. If $\varphi:\pi^{-1}(U) \To U \times G$
is a local trivialisation of $\xi$ over $U \subset N$, then let $V
= f^{-1}(U) \subset M$, and define $\psi :\pi_E^{-1}(V) \To V
\times G$ by $\psi(x,\bdelt{p}) = (x, \varphi_2(\bdelt{p}))$,
where $\varphi_2$ denotes the projection of $\varphi$ onto the $G$
factor. The condition $f(x)=\pi(\bdelt{p})$ in the definition of
the pull-back ensures that this prescription is valid.

\section{Homotopies}
Let $M$ and $N$ be topological spaces, and $f,g,h:M \To N$ be
continuous maps. We say $f$ and $g$ are \emph{homotopic} if there
exists a map $H:I \times M \To N$ so $H(0,x)=f(x)$ and
$H(1,x)=g(x)$ for all $x \in M$. The map $H$ is called a
\emph{homotopy} from $f$ to $g$. Then homotopy equivalence is in
fact an equivalence relation. Reflexivity is simple, using the
homotopy $H_f$ defined by $H_f(t,x)=f(x)$. If $H$ is a homotopy
from $f$ to $g$, define $H^{-1}$ by $H^{-1}(t,x)=H(1-t,x)$. This
appears unusual notation, in that at this point $H^{-1}$ is not an
`inverse' in any obvious sense. This will become clear later. Then
$H^{-1}$ is a homotopy from $g$ to $f$, and so homotopy
equivalence is reflexive. Transitivity is seen by defining the
composition of homotopies. If $H$ is a homotopy from $f$ to $g$,
and $K$ is a homotopy from $g$ to $h$, then we can define a
homotopy from $f$ to $h$, denoted $K \circ H$, defined by
\[
    (K \circ H)(t,x)=\left\{
        \begin{array}{ll}
            H(2t,x) & \textrm{if $t \in \left[0,\frac{1}{2}\right]$} \\
            K(2t-1,x) & \textrm{if $t \in \left[\frac{1}{2},1\right]$}
        \end{array}
        \right.
        .
\]
This is clearly continuous, and thus establishes the homotopy
equivalence of $f$ and $h$.

We would now like to specialise to the situation of smooth
manifolds and smooth maps between them. We define homotopy
equivalence the same way, requiring that the homotopy is also a
smooth map. Unfortunately the formula given above for composition
of homotopies fails, because there is no guarantee that this map
will be smooth. This is because the derivatives may not match up
when $t=\frac{1}{2}$. For our purposes, it is possible to avoid
this problem by considering a restricted class of homotopies. We
define a \emph{steady smooth homotopy} from $f$ to $g$ to be a map
$H:I \times M \To N$ so that, for some $\varepsilon > 0$,
$H(t,x)=f(x)$ for all $x \in M$ and $t \in [0,\varepsilon]$, and
$H(t,x)=g(x)$ for all $x \in M$ and $t \in [1-\varepsilon,1]$. It
is easy to see that steady smooth homotopies can be composed to
form a steady smooth homotopy. Further, if there is a smooth
homotopy from $f$ to $g$, then there is a steady smooth homotopy,
by the following prescription. Firstly define $\varphi_0:I \To I$
by
\[
    \varphi_0(s) = \left\{
        \begin{array}{ll}
        \exp\left(\frac{1}{(s-\frac{1}{3})(s-\frac{2}{3})}\right) &
        \textrm{if $s \in \left(\frac{1}{3},\frac{2}{3}\right)$} \\
        0 & \textrm{otherwise} \\
        \end{array}
        \right.
        ,
\]
and $\varphi:I \To I$ by $\varphi(t) = \frac{\int_0^t \varphi_0(s)
\mathrm{d} s}{\int_0^1 \varphi_0(s) \mathrm{d} s}$. $\varphi_0$ is
smooth, and so $\varphi$ is smooth also. Then given a smooth
homotopy $H$ from $f$ to $g$, define the steady smooth homotopy
$H'$ by $H'(t,x)=H(\varphi(t),x)$. From this result, it is seen
that we lose nothing by passing to the steady smooth homotopies,
and so homotopy equivalence is in fact an equivalence relation,
regardless of whether we use smooth homotopies or steady smooth
homotopies.

\section{Connections and paracompactness}
A paracompact manifold is a manifold with a countable basis for
its topology, or, equivalently, a manifold on which there exists a
partition of unity subordinate to any open covering. Using the
existence of partitions of unity, we can define a connection form
on any principal fibre bundle defined over the paracompact
manifold \cite{cho:amp}. Further, from this connection, we obtain
a rule for equivariant path lifting. That is, given a path in the
base space $\alpha:I \To M$, and a point $\bdelt{p} \in
\pi^{-1}(\alpha(0))$, we can form a lifted path,
$\tilde\alpha_{\bdelt{p}}:I \To P$, so $\pi \circ
\tilde\alpha_{\bdelt{p}} = \alpha$, and
$\tilde\alpha_{\bdelt{p}}(0)=\bdelt{p}$. The equivariance of the
path lifting is expressed by $\tilde\alpha_{\bdelt{p} g}(t) =
\tilde\alpha_{\bdelt{p}}(t) g$ for all $g \in G$. We can use this
path lifting to parallel transport a point of the total space.
Furthermore, the parallel transport is independent of the
parametrisation of the path, as the path lifting is defined as the
integral curve of the lifting of the tangent to the path in the
base space to the total space using the connection. That is, if
$\alpha$ is a path in $M$ and $\beta$ is a reparametrisation of
$\alpha$ such that $\alpha(0)=\beta(0)$ and
$\alpha(t_0)=\beta(t_1)$, then
$\tilde\alpha_{\bdelt{p}}(t_0)=\tilde\beta_{\bdelt{p}}(t_1)$. The
main result about parallel transport which we will rely on is that
when two paths can be composed to give a smooth path, the parallel
transports compose, in the sense that $(\widetilde{\alpha \cat
\beta})_{\bdelt{p}} = \tilde\alpha_{\tilde\beta_{\bdelt{p}}(1)}
\cat \tilde\beta_\bdelt{p}$. This fact is easily obtained from the
definition of the path lifting as an integral curve.

\section{Bundle morphisms induced by homotopies}
\subsection{A categorical approach}
In this section we prove the main result of this article, that if
$M$ and $N$ are smooth manifolds, $\xi$ is a principal fibre
bundle over $N$, and $f$ and $g$ are homotopic maps from $M$ to
$N$, then there is a principal fibre bundle isomorphism between
the pull-back bundles, $f^*\xi$ and $g^*\xi$. We'll use several
notions from category theory to organise the proof. We fix smooth
manifolds $M$ and $N$, and a principal fibre bundle $\xi$ over $N$
with structure group $G$. Firstly we consider a category
$\mathcal{H}$ whose objects are smooth maps from $M$ to $N$, and
whose morphisms are steady smooth homotopies, with composition as
defined above. We consider a second category $\mathcal{B}$, whose
objects are principal fibre bundles over $M$ with structure group
$G$, and whose morphisms are principal fibre bundle morphisms.
Notice that there is a natural map from the objects of
$\mathcal{H}$, to the objects of $\mathcal{B}$, given by taking
the pull-back of $\xi$ by the smooth map in question. This
motivates our goal of constructing a functor from $\mathcal{H}$ to
$\mathcal{B}$, that is, constructing a principal fibre bundle
morphism for each homotopy. This construction will be lead
primarily by the requirement that we obtain a functor. This
motivates the use of connections and path liftings, because of the
functorial character of the path lifting with respect to path
composition described above. The existence of functor then
establishes the final result, as follows. Since a functor takes
isomorphisms to isomorphisms, and each homotopy has an inverse
(its `time reversal'), the principal fibre bundle morphisms
obtained are also invertible, and so isomorphisms. Notice that
because of the arbitrary choice of the connection, this does not
establish a canonical isomorphism between the pull-backs.
\subsection{The construction}
Let $E$ be the total space of $f^*\xi$, and $F$ be the total space
of $g^*\xi$. From the homotopy $H$ of $f$ to $g$, we want to
construct a map $\lambda_H: E \To F$. Fix $\bdelt{e} =
(x,\bdelt{p}) \in E$. Define $\alpha : I \To N$ by $\alpha(t) =
H(t,x)$. Now according to the definition of the pull-back bundle,
$\pi(\bdelt{p}) = \alpha(0) = f(x)$. Thus, using the connection,
we can parallel transport $\bdelt{p}$ to
$\tilde\alpha_{\bdelt{p}}(1)$, and $\pi(\tilde\alpha_\bdelt{p}(1))
= \alpha(1) = g(x)$. Therefore define
$\lambda_H(\bdelt{p})=(x,\tilde\alpha_{\bdelt{p}}(1))$.
$\lambda_H$ as defined is smooth, because the solutions of
differential equations, and hence parallel transport, depend
smoothly on initial conditions. Further, $\lambda_H$ clearly
respects the fibration of $E$ and $F$, and the group action,
because $\lambda_H(\bdelt{p} g) = (x,\tilde\alpha_{\bdelt{p}
g}(1))= (x,\tilde\alpha_{\bdelt{p}}(1)g)=\lambda_H(\bdelt{p}) g$.
Thus $\lambda_H$ is in fact a principal fibre bundle morphism.
\subsection{Functoriality}
The remaining step is to ensure that the association we have given
$f \rightsquigarrow f^*\xi$, $H \rightsquigarrow \lambda_H$ is in
fact functorial, in the sense that if $H$ is a homotopy from $f$
to $g$, and $K$ is a homotopy from $g$ to $h$, so $K \circ H$ is a
homotopy from $f$ to $h$, then $\lambda_{K \circ H} = \lambda_K
\circ \lambda_H$. To check this, we fix $\bdelt{e} = (x,\bdelt{p})
\in E$, and construct $\alpha^H$, $\alpha^K$ and $\alpha^{K \circ
H}$ as
\begin{eqnarray*}
    \alpha^H(t) & = & H(t,x) \\
    \alpha^K(t) & = & K(t,x) \\
    \alpha^{K \circ H}(t) & = & (K \circ H)(t,x) \\
                       & = &
        \left\{
            \begin{array}{ll}
            H(2t,x) & \textrm{if $t \in \left[0,\frac{1}{2}\right]$} \\
            K(2t-1,x) & \textrm{if $t \in \left[\frac{1}{2},1\right]$}
        \end{array}
        \right.
\end{eqnarray*}
Thus $\alpha^{K \circ H} = \alpha^K \circ \alpha^H$, in the usual
sense of composition of paths. Then $\lambda_H(x,\bdelt{p}) =
(x,\tilde\alpha^H_{\bdelt{p}}(1))$, and
\begin{eqnarray*}
 \lambda_K(\lambda_H(x,\bdelt{p})) & = &\lambda_K(x,\tilde\alpha^H_{\bdelt{p}}(1)) \\
     & = & (x,\tilde\alpha^K_{\tilde\alpha^H_{\bdelt{p}}(1)}(1)) \\
     & = & (x,\widetilde{\alpha^K \circ \alpha^H}_{\bdelt{p}}(1)) \\
     & = & \lambda_{K \circ H}(x,\bdelt{p}).
\end{eqnarray*}

This completes the construction, and established the result. The
above result also holds for general fibre bundles.

\subsection{An application}
Our result enables a simple proof of the theorem that principal
fibre bundles over contractible base spaces are trivial. Thus, let
$\xi$ be a principal fibre bundle over $M$, a contractible smooth
manifold. Define $f,g:M \To M$ by $f(m)=m$, $g(m)=m_0$, for some
fixed $m_0 \in M$. Since $M$ is contractible, $f$ and $g$ are
homotopic, so the pull-back bundles $f^*\xi$ and $g^*\xi$ are
isomorphic. However $f^*\xi$ is simply $\xi$, and it is not hard
to see that $g^*\xi$ is trivial. This type of theorem can be used
to classify bundles. A simple extension of the theorem shows that
if the base space factors as $M=A \times B$, and $B$ is
contractible, then for every point $m \in M$, then is an open set
$U \subset A$ so $m \in U \times B$, and $\xi$ is trivial over $U
\times B$. The following is a deep theorem based on these ideas,
which is described in \cite{ish:mdgfp}. We denote by
$\mathcal{B}_G(N)$ is set of isomorphism classes of principal $G$
bundles over the base space $N$, and by $[P,Q]$ the set of
homotopy classes of maps from $P$ to $Q$.

\begin{thm}
If $G$ is a group, then there exists a \emph{universal bundle}
$\mathcal{U}G = \pfbundle{G}{PG}{}{MG}$ such that for any smooth
manifold $M$ the mapping from $[N,MG]$ to $\mathcal{B}_G(N)$
described above is one-to-one and onto. That is,
\begin{enumerate}
\item If $f,g:N \To MG$, then $f^*(\mathcal{U}G)$ and
$g^*(\mathcal{U}G)$ are isomorphic if and only if $f$ and $g$ are
homotopic maps.
\item If $\xi$ is any principal fibre bundle over $N$ with
structure group $G$ then there is some map $f:N \To MG$ so $\xi$
is isomorphic to $f^*\mathcal{U}G$.
\end{enumerate}
Different universal bundles with the group $G$ necessarily have
homotopy equivalent base spaces.\footnote{That is, they have base
spaces which are isomorphic in the category whose objects are
smooth manifolds, and whose morphisms are homotopy classes of
smooth maps.}
\end{thm}

This theorem indicates that the problem of classifying all $G$
bundles over $N$ is identical that of finding all homotopy classes
of maps from $N$ to $BG$, or, $[N,MG] \cong \mathcal{B}_G(N)$.
Particular examples include the classification result
$\mathcal{B}_{U(1)}(N) \cong [N, \Complex P^\infty] \cong
H^2(N,\mathbb Z)$, that the $U(1)$ bundles over a manifold are in
one-to-one correspondence with the second cohomology (with integer
coefficients) classes of $N$. The theory of characteristic classes
is also related to this result.

\bibliographystyle{amsplain}
\bibliography{pullbacks}
\end{document}